\documentclass[twocolumn,showkeys,preprintnumbers,amsmath,amssymb]{revtex4}

\usepackage{bm}
\usepackage{dcolumn}
\usepackage{epsfig}
\usepackage{amsfonts}
\usepackage{color}

\begin{document}


\title{{\textbf{Equation-Free Dynamic Renormalization: 
Self-Similarity in Multidimensional Particle System
Dynamics}}}

\author{Yu Zou}
 \affiliation{Department of Chemical Engineering and Program in 
Appplied and Computational Mathematics, Princeton University, Princeton, NJ 08544}
\author{Ioannis Kevrekidis}%
 \email{yannis@princeton.edu}
\affiliation{Department of Chemical Engineering and Program in 
Appplied and Computational Mathematics, Princeton University, Princeton, NJ 08544}%

\author{Roger Ghanem}
\affiliation{Department of Civil Engineering, The University of Southern 
California, Los Angeles, CA 90089}%

\date{\today}

\begin{abstract}
We present an equation-free dynamic renormalization approach
to the computational study of coarse-grained, self-similar dynamic
behavior in multidimensional particle systems.
The approach is aimed at problems for which evolution equations
for coarse-scale observables (e.g. particle density) are not explicitly 
available. 
Our illustrative example involves Brownian particles in a 2D Couette 
flow; marginal and conditional Inverse Cumulative Distribution Functions
(ICDFs) constitute the macroscopic observables of the evolving particle 
distributions.
\end{abstract}


\keywords{equation free computation, coarse dynamic 
renormalization, 
coarse time-stepper, particle system dynamics, inverse cumulative 
distribution function}

\maketitle

\section{Introduction}
Multiscale phenomena arise naturally in science and engineering. 
In many cases of current research interest, physical models are 
available
at a {\it fine, microscopic} scale (atomistic, stochastic, 
agent-based), while we want 
to study the system behavior at a {\it coarse-grained, macroscopic} 
level.
Macroscopic, coarse-grained evolution equations, even when
they {\it conceptually exist} are often 
unavailable in closed form, 
due to the lack of accurate explicit closures.
The {\it Equation-Free} computational framework 
(\cite{Theodoropoulos:00,Gear:02,Kevrekidis:03,Kevrekidis:04})
has been recently proposed for the computer-assisted study 
of precisely such complex, multiscale problems: those
for whose macroscopic behavior no {\it explicit} 
coarse-grained evolution equations are available.
Equation-free methods utilize so-called coarse time-steppers, which are 
used to numerically analyze coarse-grained behavior through 
appropriately designed short computational
experiments performed by the fine-scale models;
coarse timestepping is closely related to
optimal prediction in the work of Chorin and coworkers  \cite{Chorin:00,Chorin:02}
 (see also the discussion in \cite{Givon:04}).
Quantities necessary for traditional continuum numerical analysis 
(residuals, the action of Jacobians) are {\it estimated on demand} from these short
fine-scale computational experiments, rather than evaluated from 
closed-form macrosopic equations.
This computational ``enabling technology" has been used to perform
integration, fixed point computation, numerical stability and 
bifurcation analysis as well as control and optimization for the 
coarse-grained behavior of
a number of fine-scale model types (Lattice Boltzmann, Monte Carlo, 
Molecular Dynamics, Brownian Dynamics etc., see the references in 
\cite{Kevrekidis:03,Kevrekidis:04}).

For problems whose macroscopic behavior is characterized by {\it scale 
invariance},
dynamic renormalization methods \cite{Mclaughlin:86,LemesurierA:88}
provide tools for locating self-similar solutions
and their scaling exponents \cite{Barenblatt:96,Goldenfeld:92}.
Recently, we have combined equation-free computation with dynamic 
renormalization
to obtain coarse-grained self-similar solutions
based on direct microscopic simulation \cite{Chen:03,Chen:04}; our particular
implementation used a template-based approach 
\cite{Rowley:00,Aronson:01,Rowley:03,Siettos:03}.
Working in a frame of reference that expands (or shrinks) along with 
the
macroscopic system observables turns the self-similar problem into a
steady state one, for which fixed point techniques can be applied.
Equation free dynamic renormalization techniques are used here
to investigate macroscopic self-similarity in multidimensional particle 
system dynamics.
This is made possible through a coarse time-stepper which utilizes 
marginal and  one-dimensional conditional inverse cumulative distribution functions 
(ICDFs) as the
coarse observables of multidimensional particle systems.
Our illustrative example consists of the (self-similar) dynamics of 
Brownian particles in a Couette flow; the results are validated using
known analytical solutions.  
The paper is organized as follows: we start with a brief description of 
the coarse time-stepper in our framework. 
We then discuss 
equation-free
dynamic renormalization, present our illustrative example, demonstrate
the computation of its self-similar shapes and exponents, and conclude 
with
a brief discussion of limitations and possible extensions of the 
approach.

\section{A Coarse Time-Stepper for Multidimensional Particle Systems}

\begin{figure}
\epsfig{figure=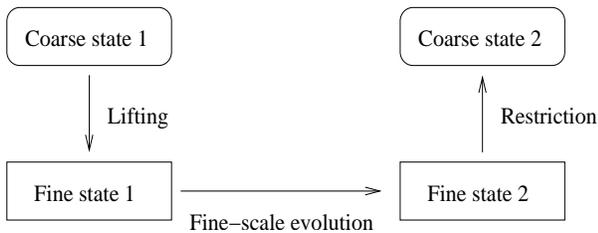,height=30mm,width=80mm,angle=0}
\caption{Schematic of the coarse time-stepper}
\label{cts:fig}
\end{figure}

In the equation-free framework, 
short, appropriately initialized computational experiments with a 
fine scale model are used to construct
the coarse time-stepper--the basic element capturing the dynamical 
interaction 
between coarse-scale observables and fine-scale model states. 
There are essentially three components in a coarse time-stepper, 
namely, 
{\it lifting}, {\it fine-scale evolution} and {\it restriction} 
\cite{Gear:02} (Fig. \ref{cts:fig}). 
Lifting is a transformation that converts coarse-scale observables to 
one or more consistent fine-scale realizations;
restriction is the transformation in the reverse direction, from fine 
scale
states to coarse observables.  
The specific manner in which these lifting and restriction
operators are implemented yields different inter-scale exchange 
protocols; one must test that the
results at the macroscopic level are relatively insensitive to small 
details
in the protocols (see the discussion in \cite{Kevrekidis:03}). 

In macroscopically multidimensional particle systems, 
particle positions are a natural component of the
fine scale state, and 
marginal and conditional ICDFs of particle positions are
good candidate coarse-scale observables.
Assuming smoothness, a finite number of marginal and conditional ICDFs 
can be used to recover the distribution of particle positions (e.g. 
through interpolation).
{\it Our coarse observables are the leading coefficients of
the projection of these ICDFs on an appropriate orthonormal
(Shifted Legendre) polynomial basis}.
Particle positions in one direction are generated using the marginal 
ICDF,
and the corresponding particle positions in the second direction are
generated using the conditional ICDFs in our lifting step.
For the restriction step, after constructing a smooth multidimensional 
CDF from particle positions by simple interpolation, marginal and conditional CDFs are first obtained, 
and then
their inverse CDFs are interpolated and their leading expansion 
coefficients computed.

\section{Equation-Free Dynamic Renormalization}

For many macroscale systems of practical interest, if the
PDEs describing their evolution are scale invariant, 
they may possess self-similar solutions \cite{Barenblatt:96,Goldenfeld:92}.
Dynamic renormalization procedures have been used for the study of such 
self-similar solutions \cite{Mclaughlin:86,LemesurierA:88,LemesurierB:88}; recently the template-based
method for studying the dynamics of problems with symmetry 
\cite{Rowley:00}
has been extended to study the dynamics of problems with scale 
invariance
\cite{Aronson:01,Siettos:03,Rowley:03}.
When scale invariant evolution PDEs are explicitly available, a 
template
based approach can be used to derive dynamical equations (termed 
``MN-dynamics")
for the renormalized self-similar solutions and similarity exponents 
\cite{Aronson:01}. 
The idea underlying this approach, especially that of employing 
template conditions,
can be used to obtain renormalized self-similar solutions 
and similarity exponents in systems 
whose macroscale governing equations are not explicitly available 
\cite{Chen:03} (see also the approach in \cite{Chorin:03}).
   
Consider a PDE of the form 
\begin{equation}
  {{\partial F} \over {\partial t}}=D_{xy}(F),
\label{eqnoriginal:eqn}
\end{equation}
where $F(x,y,t)$ is a CDF of particle positions.
We assume that for the differential operator $D_{xy}$ there exist 
constants $p$ and $a$ such that
\begin{equation}
   D_{xy}(f({x \over A},{y \over {A^p}}))=A^a D_{uv}(f(u,v)),\quad u= 
{x \over A}, \quad v={y \over {A^p}}
\label{eqnDxy:eqn}
\end{equation}
for any real function $f$, real value $A>0$ and coordinate $(x,y)$ 
(there is
no amplitude rescaling since this is a CDF).
If a  self-similar solution $F(x,y,t)$ exists, it can be written as 
\begin{equation}
   F(x,y,t)=U({x \over (cs)^\alpha}, {y \over (cs)^{\alpha p}};c),
\label{eqn7:eqn}
\end{equation}
where $c$ is a constant parametrizing the family of 
renormalized shapes, $s=t-t_0$, ($t_0$ is the blowup time for 
problems with finite time singularities, whether forward or backward in time).

Therefore, 
\begin{equation}
  \alpha a = -1,
\label{eqnalphaa:eqn}
\end{equation}
and $U$ satisfies the PDE,
\begin{equation}
  - \alpha u U_u- \alpha p v U_v = c^{-1} D_{uv}(U),
\label{eqnselfsimilar:eqn}
\end{equation}
where $u=x/(cs)^\alpha$, $v=y/(cs)^{\alpha p}$. For $D_{xy}$ satisfying 
Eqn.(\ref{eqnDxy:eqn}), 
the constant $a$ is determined by $D_{xy}$ itself and the similarity 
exponent $\alpha$ can be obtained by Eqn.(\ref{eqnalphaa:eqn}).

If the equation is not explicitly available, one cannot analytically 
obtain
the exponents $a$ and $p$; tests have to be devised for finding these
constants --and thus testing the scale invariance of the operator-- 
before
we embark upon the computation of the self-similar solutions 
themselves.

For $D_{xy}$ satisfying the (unknown) equation (\ref{eqnDxy:eqn}), 
the constants $p$ and $a$ can be obtained {\it using a black box simulator} 
of the
equation as follows:
Since the unknown Eqn. (\ref{eqnDxy:eqn}) is valid for {\it any} real 
function $f$, 
let $f$ be a test function, and let us choose the two points $(u_1,v_1)$ 
and $(u_2,v_2)$
and let
$(x_1,y_1)=(u_1 A,v_1 A^p)$ and $(x_2,y_2)=(u_2 A,v_2 A^p)$, 
where $A$ is arbitrarily chosen as a positive real value.
Then the following two relations should hold for essentially arbitrary
test functions $f(x,y)$ and points $(u_i,v_i), i=1,2$:
\begin{eqnarray}
   D_{xy}(f({x \over A},{y \over {A^p}}))(x_1,y_1) &=& A^a 
D_{uv}(f(u,v))(u_1,v_1),  \nonumber \\
   D_{xy}(f({x \over A},{y \over {A^p}}))(x_2,y_2) &=& A^a 
D_{uv}(f(u,v))(u_2,v_2).
\end{eqnarray}
Comparing the above two equations, we have
\begin{equation}
   {{D_{xy}(f({x \over A},{y \over {A^p}}))(x_1,y_1)} \over 
{D_{xy}(f({x \over A},{y \over {A^p}}))(x_2,y_2)}}= {{D_{uv}(f(u,v))(u_1,v_1)} 
\over {D_{uv}(f(u,v))(u_2,v_2)}}.
\label{eqnconp:eqn}
\end{equation}
Eqn.(\ref{eqnconp:eqn}) is solved (using Newton's method) for the 
constant $p$.
When $D_{xy}$ is not explicitly available, $D_{xy}(f(x,y))$ can be 
obtained
by running the micro-simulator for a short time and 
numerically estimating the derivative ${\partial f} \over {\partial 
t}$. 
Given $p$, the constant $a$ is calculated by 
\begin{equation}
    a=\log_{A}{{D_{xy}(f({x \over A},{y \over {A^p}}))(x_1,y_1)} \over 
{{D_{uv}(f(u,v))(u_1,v_1)}}}. 
\label{eqnasolve:eqn}
\end{equation}
Clearly, other test functions, as well as conditions evaluated at other 
points 
(or nonlocal versions of the conditions using templates) can be used;
care must be taken also to ensure the finiteness of the estimated 
quantities.

Given $a$ and $p$, to determine the self-similar shape of the solution, 
we consider the general scaling
\begin{equation}
  F(x,y,t)=\omega ( {x \over A(t)}, {y \over {A(t)^p}}, t),
\label{eqnrescaling:eqn}
\end{equation}
where $A(t)$ is an unknown function. The PDE becomes
\begin{equation}
  \omega_t-{ A_t \over A} u \omega_u - {{p A_t} \over A} v \omega_v = 
A^a D_{uv}(\omega).
\label{eqnomega:eqn}
\end{equation}
Evidently, $U$ and $\omega$ are both rescaled CDFs.

A single template condition coupled with equation (\ref{eqnomega:eqn})
is needed to solve for both $\omega(u,v,t)$ and  $A(t)$
at every time step. 
In this paper we choose this template condition to be
\begin{equation}
  \omega(e,\infty,t)=m, \quad e<0, \quad  0<m<0.5,
\end{equation}
where $e$ and $m$ are both (essentially arbitrary) constants. 
This template condition has the following physical meaning: the 
u-coordinates corresponding to the rescaled marginal CDF $\omega_U=m$ 
have the same value $e$ for all $t$.

Applying this template to Eqn.(\ref{eqnomega:eqn}) 
and assuming ${{\partial \omega} \over {\partial v}}(e,v,t)$ 
decays exponentially as $v \rightarrow \infty$, we have
\begin{equation}
  A_t e {{\partial \omega} \over {\partial u}} (e,\infty,t) + A^{a+1} 
D_{uv}(\omega)(e,\infty,t) = 0.
\label{eqncouple:eqn}
\end{equation}  

Equations (\ref{eqnomega:eqn}) and (\ref{eqncouple:eqn}) are solved for 
the rescaled CDF $\omega(u,v,t)$ and rescaling variable $A(t)$ if the operator 
$D_{xy}$ is explicitly known.
As the time $t \rightarrow \infty$, $\omega$ may approach a steady 
state, 
which is a stable self-similar shape for the solutions to 
Eqn.(\ref{eqnoriginal:eqn}). 

The value for $\alpha$ can be calculated in the long-time limit (i.e.,
after $\omega$ reaches steady state). 
Indeed, let $t_1$ and $t_2$ be distinct times after $\omega$ reaches 
the steady state, then

\begin{equation}
   \alpha= { {t_2-t_1} \over { { A(t_2) \over A_t(t_2)} -{ A(t_1) \over 
A_t(t_1)}} }.
\label{alphacal:eqn}
\end{equation}

\begin{figure}
\epsfig{figure=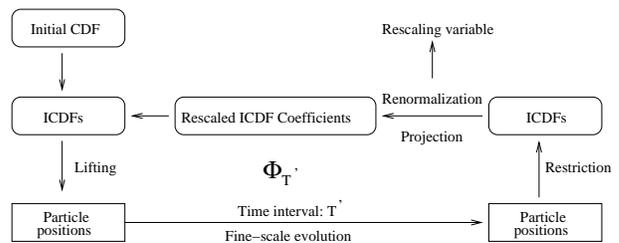,width=0.45\textwidth}
\caption{Illustration for the coarse renormalization}
\label{crm:fig}
\end{figure}

In many cases the macroscale equation for the CDF of particle positions 
may not be explicitly known. 
However, the template conditions can still be used to renormalize the 
CDF itself, evolved 
via microscale models, and rescaling variables are obtained during the 
course of renormalization.
In these cases, a coarse time-stepper 
is employed to evolve the coarse system using the 
the marginal and conditional ICDFs (actually, the leading coefficients
of their projection on an appropriate basis) as coarse observables.
The procedure for equation-free renormalization, depicted in Fig. 
\ref{crm:fig},
consists of the following six steps:
  
\begin{enumerate}
\item Generate the marginal and conditional ICDFs according 
to the initial CDF or according to the coefficients of dominant modes 
of the ICDFs.
\item Generate particle positions based on these ICDFs using the {\it 
lifting} procedure in the coarse time-stepper.
\item Evolve particle positions over a time interval $T'$ using the 
fine-scale model.
\item Obtain ICDFs from particle positions using the {\it restriction} 
procedure in the coarse time-stepper.
\item Rescale the marginal ICDF according to the template condition and 
obtain the rescaling variable $A$. 
We then rescale the conditional ICDFs by $A^p$ (recall that $p$ has been 
independently computed through Eqn. (\ref{eqnconp:eqn})). 
This step can be justified by Eqn.(\ref{eqnrescaling:eqn}). 
Indeed, obtaining the rescaled solution $\omega_{k+1}$ from $\omega_k$ 
via the dynamics (\ref{eqnomega:eqn}) and (\ref{eqncouple:eqn}) is 
equivalent to starting 
from the initial condition $\omega_k$, running the original dynamics 
for a while to get $F_{k+1}$, 
obtaining the rescaling variable $A$, and rescaling $F_{k+1}$ by $A$ 
and $A^p$ respectively in the 
$x$ and $y$ directions. 
\item Project the rescaled ICDFs onto the appropriate orthonormal 
basis and obtain the coefficients of its dominant modes. Go back to step 1.
\end{enumerate} 

The above procedure can be viewed as an iterative algorithm to 
solve for the fixed point of a nonlinear operator $\Phi_{T'}$,
written as,

\begin{equation}
{\bm \alpha}-\Phi_{T'} ({\bm \alpha}) = {\bm 0}\ .
\label{eqn8:eqn}
\end{equation}
This fixed point can be written in component form as,
$\alpha^r_{i,p},i=1,\cdots,M+1,p=0,\cdots,P$ \cite{Chen:03}, or

\begin{eqnarray}
{\bm \alpha} &=& 
(\alpha^r_{1,0},\alpha^r_{1,1},\cdots,\alpha^r_{1,P},\alpha^r_{2,0},\alpha^r_{2,1},\cdots,\alpha^r_{2,P}, \nonumber \\
&\cdots&,\alpha^r_{M+1,0},\alpha^r_{M+1,1},\cdots,\alpha^r_{M+1,P})^T,  
\quad \nonumber
\end{eqnarray}
where $\alpha^r_{1,j},j=0,1,\cdots,P$ stands for the projection of the 
marginal ICDF 
onto a $j$th-order mode, $\alpha^r_{i,j},i=2,\cdots,M+1;j=0,1,\cdots,P$ 
the projection 
of the $i-1$th conditional ICDF onto a $j$th-order mode, and the 
superscript $r$ refers 
to the fact that this is in fact the renormalized self-similar shape. 
These coefficients correspond to the renormalized ICDFs and CDF of the 
multidimensional particle system. 
Equation (\ref{eqn8:eqn}) can also be solved using any numerical algorithm 
such as 
direct iteration or matrix-free implementations of
Newton's method (e.g. Newton-GMRES, \cite{Kelley:95}).

\section{Brownian Particles in Couette Flow}
We now illustrate our computational approach to self-similarity in a two-dimensional Brownian model of 
particle dispersion
in Couette flow \cite{Panton:96}.
Let $X(t)$ and $Y(t)$ represent particle positions 
in the $x$ and $y$ directions respectively at time $t$ on the plane. 
The particle positions evolve according to

\begin{equation}
  dX(t)=Dd\omega_X (t), \quad dY(t)=Xdt,
\label{eqn9:eqn}
\end{equation}
where $\omega_X (t)$ is a Wiener process \cite{Gihman:72} and $D$ is 
the diffusion coefficient. The discretized dynamics of (\ref{eqn9:eqn}) 
is given by \cite{Milshtein:74}

\begin{equation}
   \quad X_{k+1}=X_k+D\eta_{X,k}\sqrt{\Delta t}, \quad 
Y_{k+1}=Y_k+X_kdt,
\label{eqn10:eqn}
\end{equation}
where $\eta_{X,k}$ are i.i.d. standard Gaussian random variables. 

The dynamics (\ref{eqn9:eqn}) represent the motion of 
particles {\it which only diffuse in the $x$-direction} in a Couette flow
(see the Discussion for both $x$ and $y$ diffusion).
It can be shown that the coarse-scale dynamics for the PDF, 
$P_{XY}(x,y,t)$, of a particle position, corresponding to the fine-scale dynamics 
(\ref{eqn9:eqn}), is governed by the following equation \cite{Majda:99}

\begin{equation}
  {{\partial P_{XY}} \over {\partial t}} + x {{\partial P_{XY}} \over 
{\partial y}} = {{D^2}\over 2} {{\partial^2 P_{XY}} \over {\partial 
x^2}},
\label{eqn11:eqn}
\end{equation} 
where $P_{XY}$ is assumed to be 2nd-order differentiable. Hence the 
dynamics for the CDF, $F_{XY}(x,y,t)$, associated with (\ref{eqn11:eqn}) 
is given by

\begin{equation}
  {{\partial F_{XY}} \over {\partial t}} + x {{\partial F_{XY}} \over 
{\partial y}} - \int_{-\infty}^x {{\partial F_{XY}} \over {\partial y}} 
dx = {{D^2}\over 2} {{\partial^2 F_{XY}} \over {\partial x^2}} \ .
\label{eqn12:eqn}
\end{equation}

In the above equation, the operator $D_{xy}$ is written as
\begin{equation}
  D_{xy}= - x {{\partial } \over {\partial y}} + \int_{-\infty}^x 
{{\partial} \over {\partial y}} dx + {{D^2}\over 2} {{\partial^2 } \over 
{\partial x^2}}.
\end{equation}

This operator satisfies the invariance property (\ref{eqnDxy:eqn}), 
with constants $p=3$ and $a=-2$. 

The analytical self-similar solution to the PDF equation 
(\ref{eqn11:eqn}) is
\begin{equation}
  P_{XY}(x,y,t)={\sqrt{3} \over {\pi {D^2} (t-t_0)^2 }} e^{- \left( 
{{6(y - 0.5x (t-t_0))^2} \over {D^2 (t-t_0)^3}} + {{x^2} \over {2D^2 
(t-t_0)}} \right)}, 
\label{eqn102:eqn}
\end{equation} 
where $t_0$ is the blowup time (forward or backward in time); the corresponding 
CDF self-similar solution 
to (\ref{eqn12:eqn}) is
$$
  F_{XY}(x,y,t)=
$$
\begin{equation}
{\sqrt{3} \over {\pi {D^2} (t-t_0)^2}} \int_{-\infty}^x \int_{-\infty}^y e^{- ( {{6(y - 
0.5x(t-t_0))^2} \over {D^2 (t-t_0)^3}} + {{x^2} \over {2D^2 (t-t_0)}} )} dy dx.
\label{eqn103:eqn}
\end{equation} 
Let $u'={x \over {(c(t-t_0))^{1/2}}}$ and $v'={y \over 
{(c(t-t_0))^{3/2}}}$; then 
$$
 F_{XY}(x,y,t) = F_{UV}(u',v')= \nonumber  
$$
\begin{equation}
 { {\sqrt{3}} \over {\pi D^2/c^2}} \int_{-\infty}^{u'} \int_{-\infty}^{v'} e^{- ( {{6(v - 
0.5u/c)^2} \over {D^2/c^3}} + {{u^2} \over {2D^2/c}} )} dv du.
\label{eqn104:eqn}
\end{equation}
Hence for the integro-differential equation (\ref{eqn12:eqn}), 
the similarity exponent $\alpha$ in (\ref{eqn7:eqn}) is $\alpha=1/2$. 
For the CDF in (\ref{eqn104:eqn}), its std.'s in two directions 
and correlation are $\sigma_X=D/c^{1/2}$, 
$\sigma_Y=D/(\sqrt{3}c^{3/2})$, 
and $\rho_{XY}=\sqrt{3}/2$, respectively. 

\section{Numerical Results} 
We now compute the self-similar solutions {\it without} the macroscopic 
equations, based only on the particle simulator.
The invariance property of the (presumed unavailable) macroscale 
differential operator $D_{xy}$ 
has to be established first. 
Then the fixed-point algorithm is used to solve for the renormalized 
CDF steady-state shape,
and the similarity exponent is computed.
The diffusion coefficient $D$ and simulation time step $\Delta t$ 
in the fine-scale model are set to $5.0cm/s^{1/2}$ and $0.01s$, 
respectively. 
An ensemble of 9000 ($N=9000$) particles is typically used in the fine-scale 
simulations. 

Without the macroscale equation (\ref{eqn12:eqn}), 
we use Newton's method to solve equation (\ref{eqnconp:eqn}) for 
$p$. 
Our test function $f$ is a 2-dimensional joint Gaussian, 
$f(x,y)=1/4.5^2N(x/4.5)N(y/4.5)$, where $N(x)$ and $N(y)$ are standard 
Gaussian distributions. 
The positive real number $A=2.0$. 
The two coordinates $(u_1,v_1)$ and $(u_2,v_2)$ in 
(\ref{eqnconp:eqn}) 
are chosen as $(-2.5,-2.5)$ and $(3.5,3.5)$, respectively. 
To reduce fluctuations of values for the operator $D_{xy}$, 
200 copies of values for $D_{xy}$ and $D_{uv}$ are averaged in the 
computation.

Starting from the initial value $p_0=6.0$, iterative values for $p$ are 
stabilized at $3.0$ after $3$ iterations. 
Accordingly, the converged value for $a$ is approximately $-2.0$ 
(Table \ref{tabpa1:tab})
(other test functions and test point selections gave similar results).

\begin{table}
\begin{tabular*}{0.45\textwidth}%
   {@{\extracolsep{\fill}}cccr}
   No. of Iterations &    $p$   &    $a$ \cr
   \hline
   0                  &   6.0  & -3.86287  \cr
   1                  &  2.79441  & -1.83047  \cr
   2                  &  2.97849  & -1.97927  \cr
   3                  &  3.00132  & -2.00114  \cr
   4                  &  2.99804  & -1.99732  \cr
   5                  &  2.99670  & -1.99612  \cr
   6                  &  2.99441  & -1.99675  \cr
   7                  &  2.99343  & -1.99113  \cr
   8                  &  2.99400  & -1.99459  \cr
   \hline
\end{tabular*}
\caption{Iterations in the computation of the constants $p$ and $a$}
\label{tabpa1:tab}
\end{table}

This suggests that the unavailable differential operator $D_{xy}$ 
corresponding to the 
microsimulator (\ref{eqn9:eqn}) does indeed possess the invariance property 
(\ref{eqnDxy:eqn})
for constants $p=3.0$ and $a=-2.0$.

To determine the self-similar shape, the template condition for the $x$ direction 
was chosen 
as $\omega(-2.266,\infty,t)=0.4$, i.e., the x-coordinates 
corresponding 
to the renormalized marginal CDF $\omega_U=0.4$ always has the same 
value, $-2.266$cm. 
For the analytical solution (\ref{eqn104:eqn}) and our template
$c=0.3125sec^{-1}$ and the corresponding std.'s 
are $\sigma_X=4\sqrt{5}$cm and $\sigma_Y=12.8\sqrt{15}/3$cm, 
respectively.  

\begin{figure}
\begin{minipage}{4cm}
\epsfig{figure=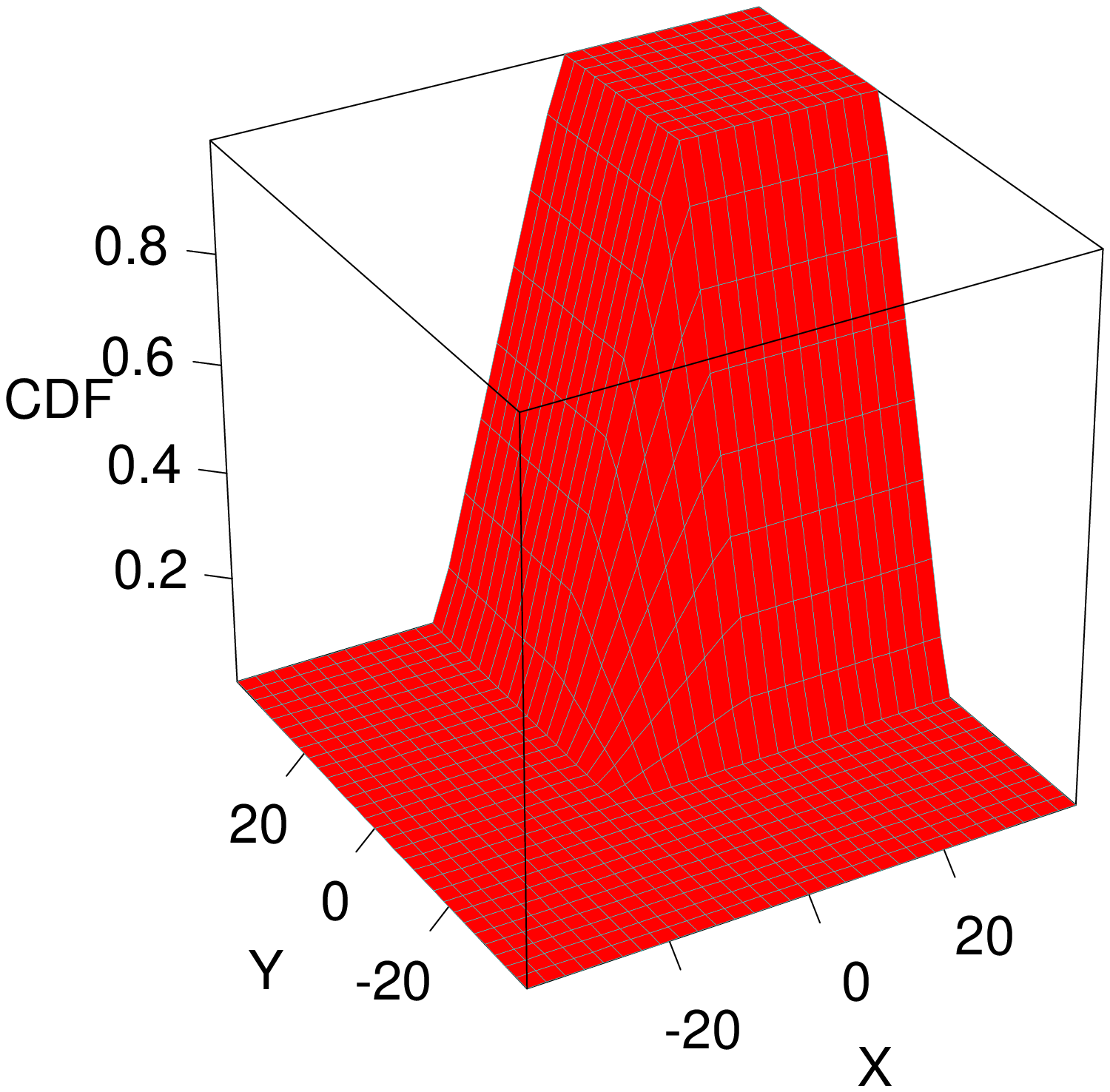,width=\textwidth}
\end{minipage}
\hfill
\begin{minipage}{4cm}
\epsfig{figure=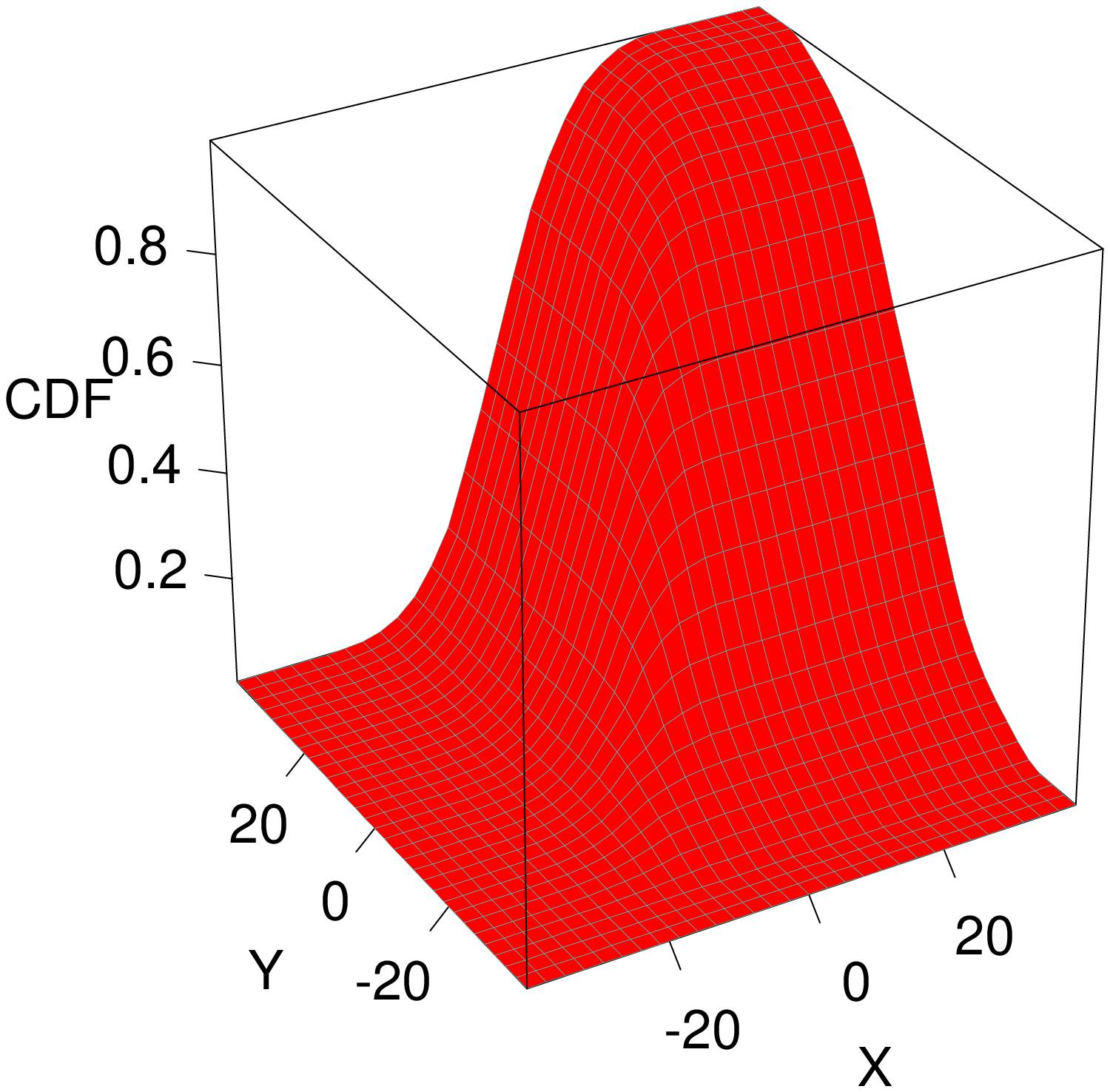,width=\textwidth}
\end{minipage}
\hfill
\begin{minipage}{4cm}
\epsfig{figure=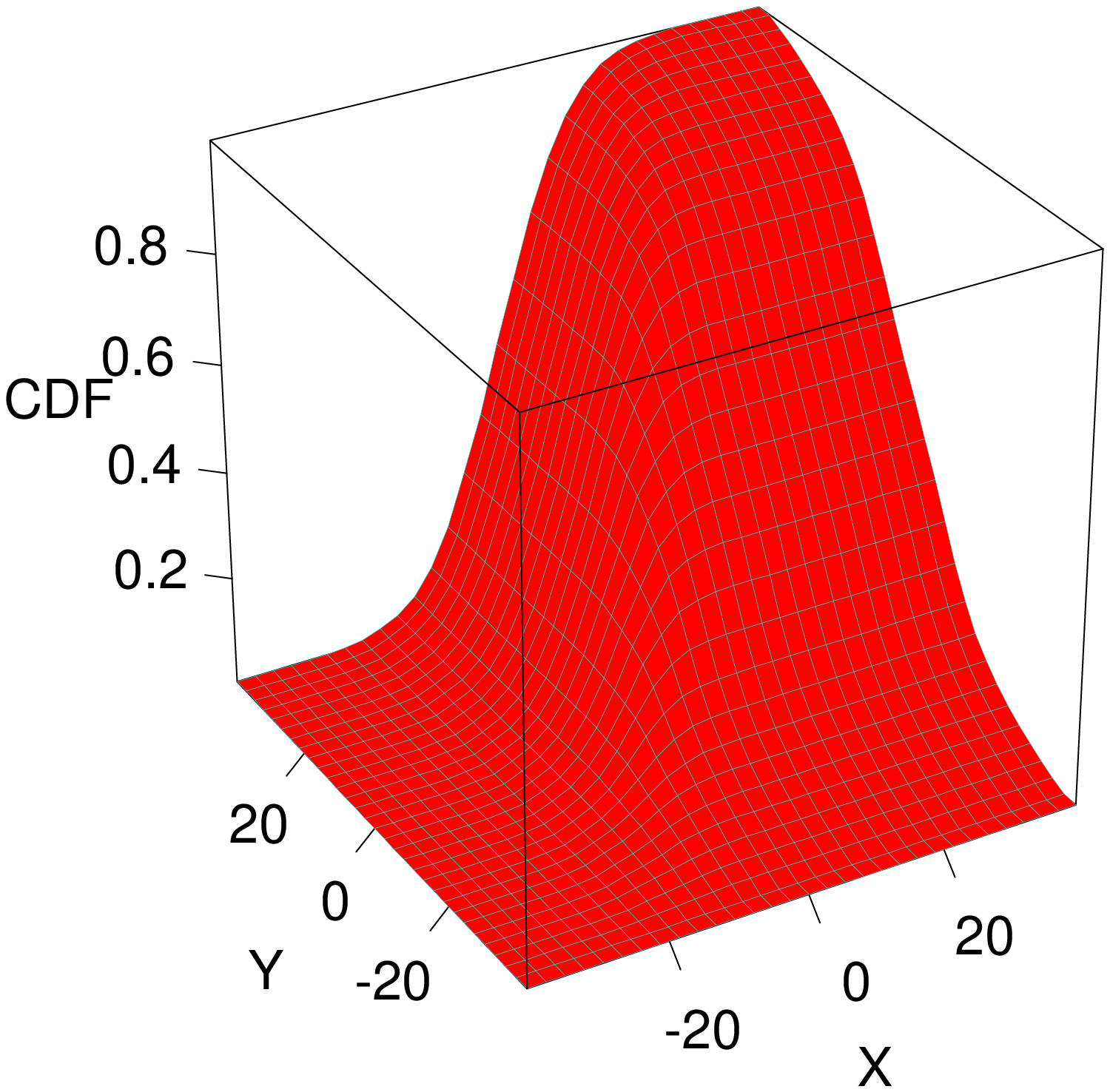,width=\textwidth}
\end{minipage}
\hfill
\begin{minipage}{4cm}
\epsfig{figure=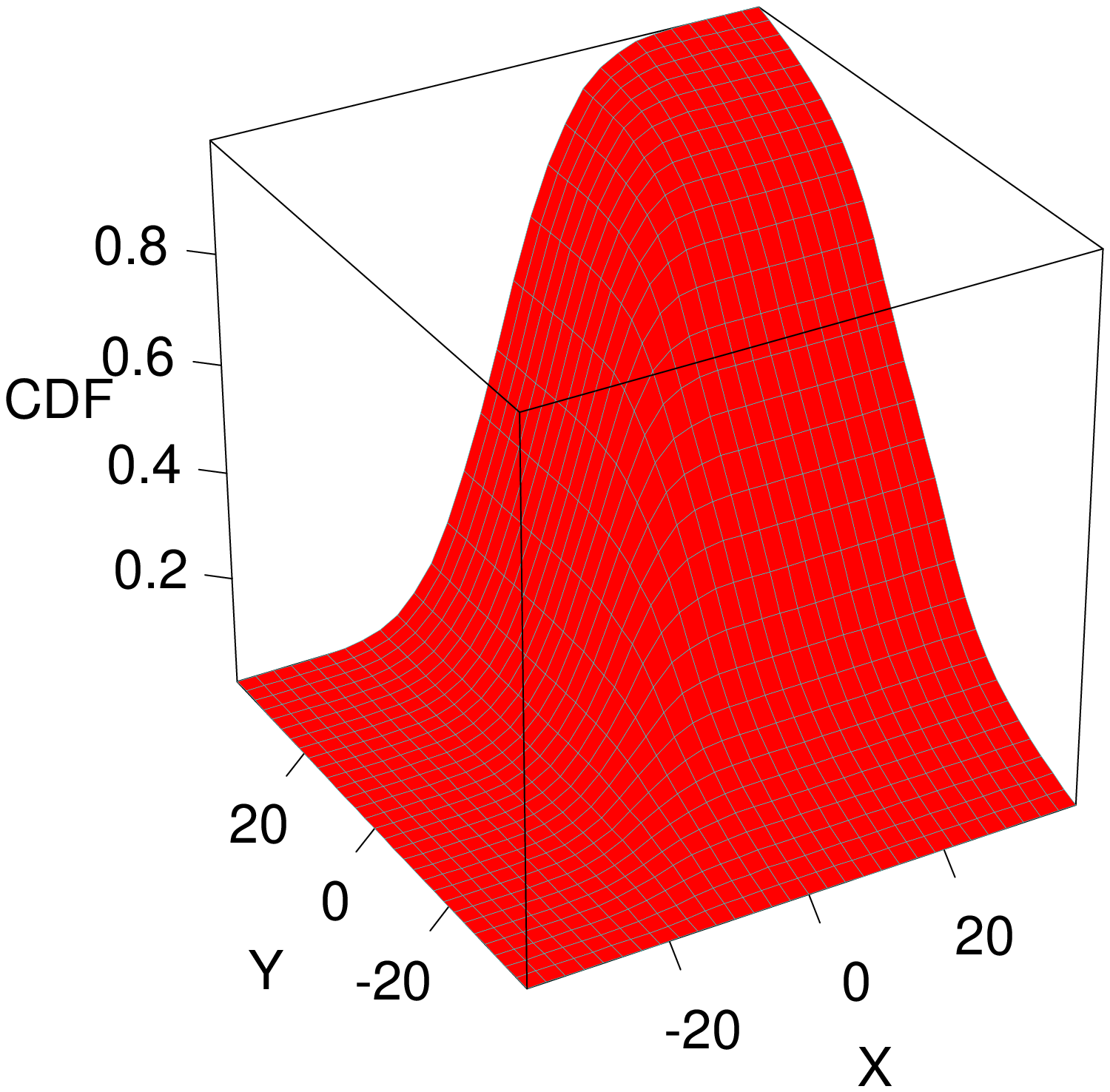,width=\textwidth}
\end{minipage}
\hfill
\caption{Renormalized CDFs ; top left: initial CDF, top right: 
renormalized CDF after 2nd iteration, bottom left: renormalized CDF after 4th 
iteration, bottom right: renormalized CDF after 6th iteration (see text).} 
\label{fig8:fig}
\end{figure}

A uniform particle CDF over the space domain 
$(-10cm,10cm)\times(-10cm,10cm)$ is used as the initial condition
for our equation-free fixed point algorithm.
Direct iteration is used to solve the fixed point of equation 
(\ref{eqn8:eqn}). 
The time interval $T'$ is $150 \Delta t$. 
The number of conditional ICDFs is 20 ($M=20$) and the basis for the 
both types of ICDF is the 
shifted Legendre polynomials of order up to and including 5 ($P=5$). 
In this simulation, 200 copies of 5000 particle positions are 
generated according to the ICDF mode
coefficients at the beginning of each iteration and let to evolve for 
time $T'$. 
The mode coefficients at the end of each iteration are obtained by 
averaging over these 200 copies. 
After the 2nd, 4th and 6th iterations, mode coefficients of 
renormalized ICDFs are used to generate particle positions, 
out of which the CDFs are computed and plotted respectively in Fig. 
\ref{fig8:fig}. 
It can be seen that the renormalized solutions reach a steady state 
quickly. 

\begin{figure}
\epsfig{figure=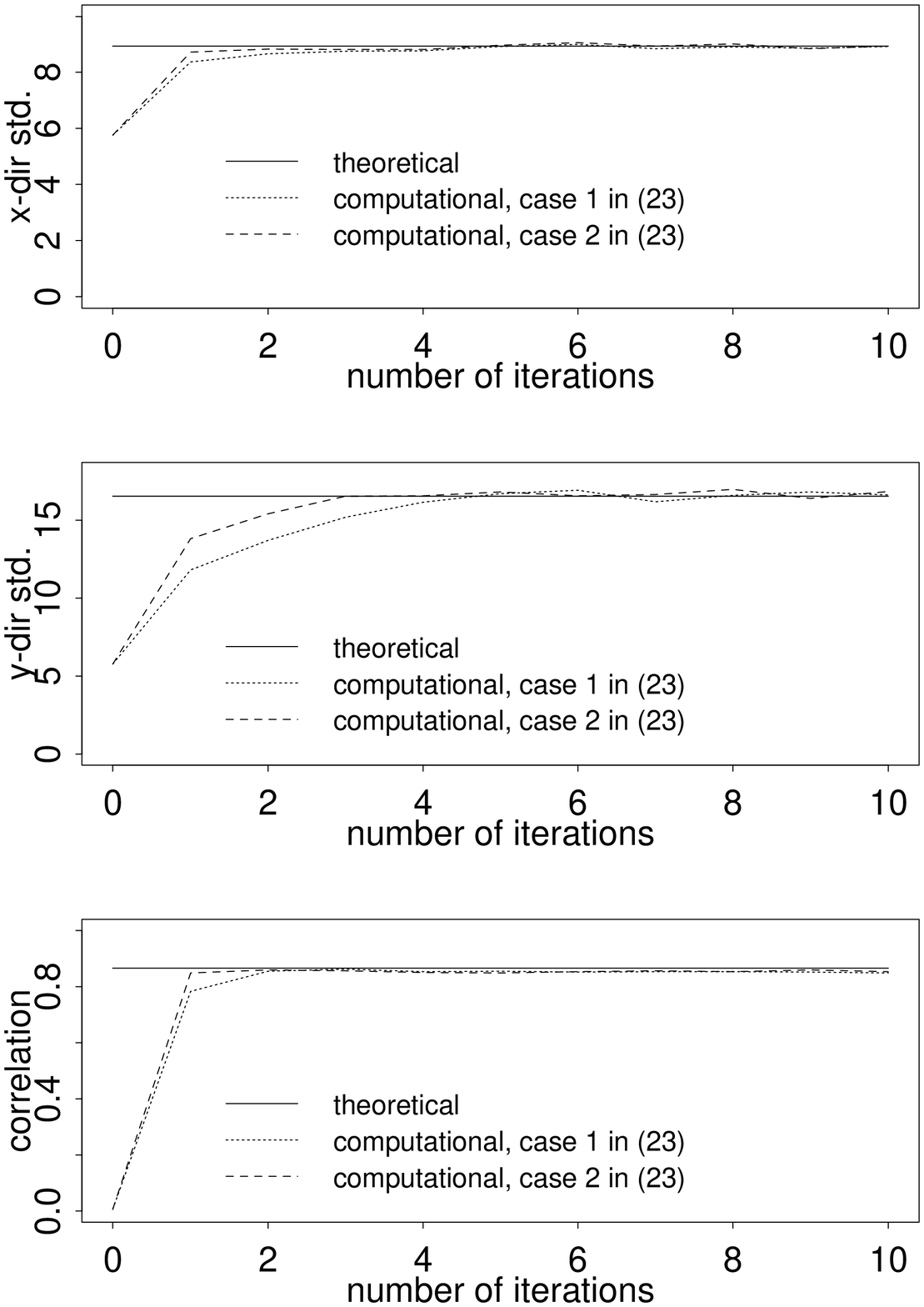,width=0.35\textwidth}
\caption{Computational vs. theoretical values for the standard deviations 
and the correlation of the self-similar shapes. Case 1, 2 in (23) (see text).}
\label{fig10:fig}
\end{figure}

To verify the self-similar shape of the solution, 
the std.'s and correlation of the simulated self-similar shape 
are compared with those of the analytical self-similar shape. 
The std.'s and correlations of the rescaled CDFs are calculated 
via the ensemble particle positions corresponding to these CDFs. 
The comparison is shown in Fig. \ref{fig10:fig}. 
The std.'s and correlations of rescaled CDFs approach those
of the analytical self-similar shape, which means that the rescaled CDF 
quickly approaches a member of the family of theoretical 
self-similar shapes 
expressed by Equation (\ref{eqn104:eqn}).
 
As the renormalized CDF $\omega(u,v,t)$ reaches the steady state, 
we can set this CDF as the initial condition and run the microscale 
dynamics 
(\ref{eqn10:eqn}) for two more loops with $t_1=150\Delta t$ and 
$t_2=300 \Delta t$. 
The rescaling variable $A(t)$ is listed in Table \ref{tabA:tab}. 
Note that $A(t)=1$ at $t=0$. By Equation (\ref{alphacal:eqn}), 
the similarity exponent $\alpha$ is approximately estimated as $0.461$, 
which is within $8\%$ of the theoretical value $1/2$.

\begin{table}
\begin{tabular*}{0.45\textwidth}%
   {@{\extracolsep{\fill}}cccr}
   $t(sec)$           &   $A(t)$  &  $A_t(t)$ \cr
   \hline
   0                  &  1.00000  &    -      \cr
   1.5                  &  1.20872  &  0.13915  \cr
   3.0                  &  1.38243  &  0.11581  \cr
   \hline
\end{tabular*}
\caption{Values of the rescaling variable $A(t)$ in Eqn.(\ref{alphacal:eqn})}
\label{tabA:tab}
\end{table}

We now study the effect of varying templates and reporting horizons
in the fixed-point operator $\Phi_{T'}$, on our results;
different templates should give different members in
the family of Eqn. (\ref{eqnselfsimilar:eqn}), while $T'$ should not
affect the computed fixed points.

Four combinations of template condition and evolution times including the one 
above are investigated:

\begin{eqnarray}
1.& \omega(-2.266,\infty,t)=0.4,& T'=150\Delta t;  \nonumber
                                          \\
2.& \omega(-2.266,\infty,t)=0.4,& T'=250\Delta t;  \nonumber
                                          \\
3.& \omega(-0.227,\infty,t)=0.4,& T'=150\Delta t;  
                                          \\
4.& \omega(-0.227,\infty,t)=0.4,& T'=250\Delta t.  \nonumber 
\end{eqnarray} 

\begin{figure}
\epsfig{figure=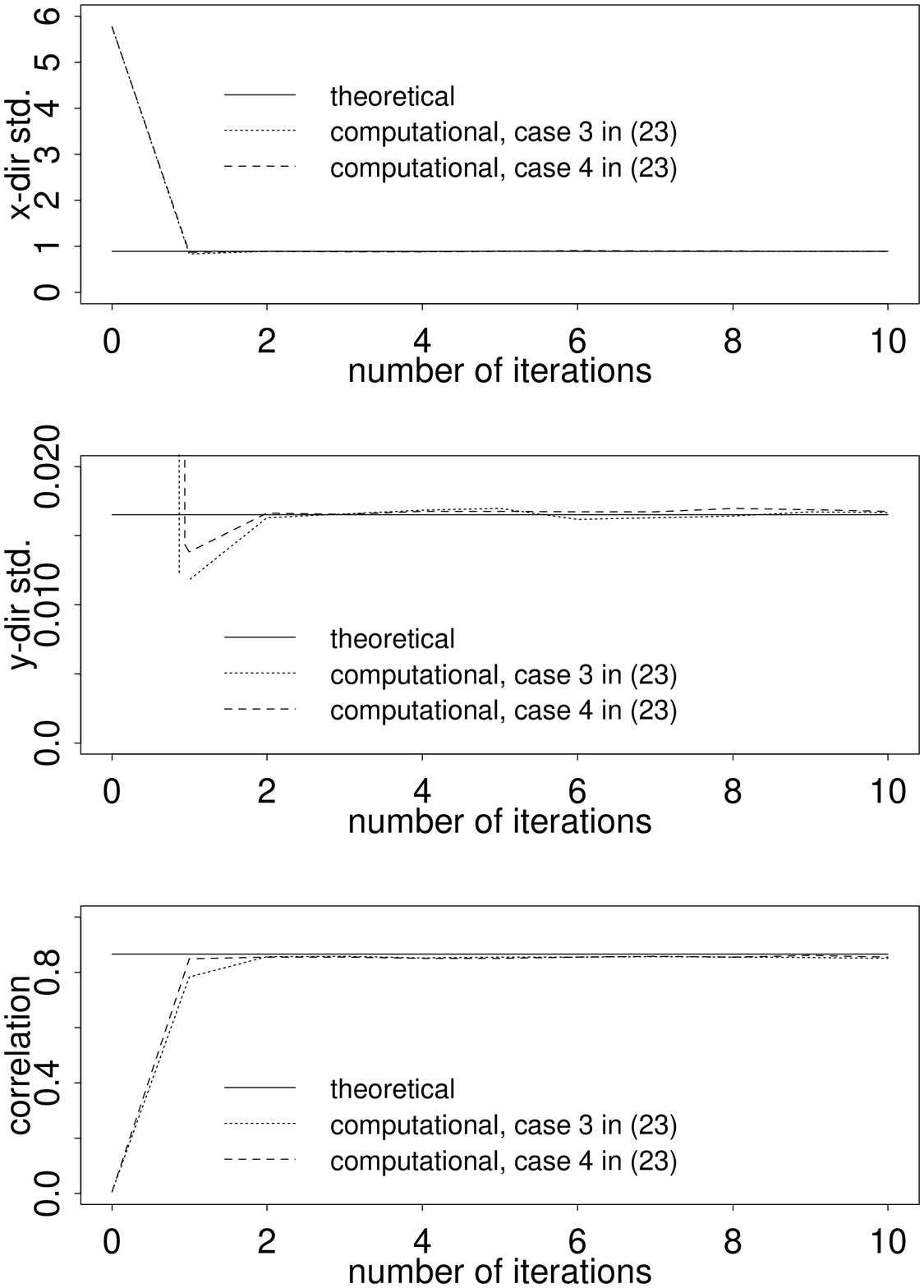,width=0.35\textwidth}
\caption{Computational vs. theoretical values for the standard deviations 
and the correlation of the self-similar shapes. Case 3, 4 in (23) (see text).}
\label{selfsimtemevocom2:fig}
\end{figure}

The iterative values of std.'s and correlation for 
the four cases are shown in figures \ref{fig10:fig} and 
\ref{selfsimtemevocom2:fig}. 
Comparison with theoretical calculations shows that variation of 
templates and evolution times indeed does not cause deviation of 
the converged rescaled CDF from the family of self-similar solutions.    
     
\section{Conclusions}

We presented a coarse dynamic renormalization technique, based on a 
coarse time-stepper using marginal and conditional ICDFs as 
coarse-scale observables, for the computer-assisted analysis of multidimensional 
self-similar particle systems. 
For coarse-grained differential operators that possess a 
scale-invariance property of the type arising in
our numerical example, a single template condition is
necessary for rescaling the ICDFs. 
In our example, We confirmed that the steady-state 
shape of the renormalized CDFs of particle positions is not affected by the 
template condition and evolution time interval in our coarse time-stepper. 

The equation-free dynamic renormalization technique was applied
only to a simple 2D Brownian particle system in this paper. 
In our computations we had already factored out translational 
invariance; we also knew the particular coordinate system
(the $x$ and $y$ axis) in which the macroscopic differential
operator is scale invariant. Finding this information, if
unknown, can be incorporated in the test procedure for
scale invariance.
For higher-dimensional diffusive particle systems, this technique may still be used.
In these cases, conditional ICDFs in third or higher dimensions have to be
utilized, thus rendering the computations slightly more complicated, yet still qualitatively 
similar to the ones presented here.
Using appropriately selected template conditions,
this technique can also be used to approximate coarse-grained
{\it asymptotically} self-similar shapes;
this is the case for our Couette example when the particles diffuse in {\it both} the
$x$ and $y$ directions (work in progress).
Several other examples of coarse self-similarity
(e.g. models of glassy dynamics, KPZ-type evolution of interfaces,
dynamics of energy spectra in randomly forced PDEs) are also being 
explored.

\clearpage

\bibliographystyle{plain}
\bibliography{PREver1}

\end{document}